\documentclass[11pt]{amsart}
\usepackage{amsmath}
\usepackage{amssymb}
\usepackage[a4paper]{geometry}
\newtheorem{teo}{Theorem}[section]
\newtheorem{pro}{Proposition}[section]

\newtheorem{cor}{Corollary}[section]

\theoremstyle{remark}
\newtheorem*{nota}{Remark}

\newcommand{\re}{\mathbb{R}}

\newcommand{\g}{\Gamma}
\newcommand{\ja}{P_n^{\alpha,\beta}}
\usepackage{mathrsfs}
\title{Free Martingale polynomials For Stationary Jacobi Processes} 
\date{\today}
\begin{document}
\maketitle
\centerline{N. DEMNI \footnote{Laboratoire de Probabilit\'es et Mod\`eles Al\'eatoires, Universit\'e de Paris VI, 4 Place Jussieu, Case 188, F-75252 Paris Cedex 05, e-mail: demni@ccr.jussieu.fr.\\
keywords : stationary free Jacobi process, multiplicative renormalization method, Tchebycheff polynomials.}}
\begin{abstract}
We generalize a previous result concerning free martingale polynomials for the stationary free Jacobi process of parameters $\lambda \in ]0.1], \theta = 1/2$. Hopelessly, apart from the case $\lambda = 1$, the polynomials we derive are no longer orthogonal with respect to the spectral measure. As a matter of fact, we use the multiplicative renormalization method to write down its corresponding orthogonal polynomials as well as the orthogonality measure associated with the martingale polynomials. We finally give a realization of the spectral measure of the free stationary Jacobi process by means of the corresponding one mode interacting Fock space.  
\end{abstract}
\section{Preliminaries}
Let $(\mathscr{A},\phi)$ a $W^{\star}$-non commutative probability space. Easily speaking, $\mathscr{A}$ is a unital von Neumann algebra and $\phi$ is a tracial faithful linear functional (state). In a previous work (\cite{Dem}), we defined, via matrix theory, and studied a two parameters-dependent self-adjoint free process, called free Jacobi process. Our focus will be on a particular case called the stationary Jacobi process since its spectral distribution does not depend on time. It is defined as $J_t: = PUY_tQY_t^{\star}U^{\star}P$ where 
\begin{itemize}
  \item $(Y_t)_{t \geq 0}$ is a free multiplicative Brownian motion (see \cite{Bia}). 
  \item $U$ is a Haar unitary operator in $(\mathscr{A}, \Phi)$.
  \item $P$ is a projection with $\Phi(P) = \lambda \theta \leq 1$, \, $\theta \in ]0, 1]$.
  \item $Q$ is a projection with $\Phi(Q) = \theta$.
  \item $QP = PQ = \left\{\begin{array}{ccc}P & \textrm{if} \quad \lambda \leq 1 \\  Q & \textrm{if} \quad \lambda > 1\end{array}\right.$
  \item $\{U, U^{\star}\}$ and $\{P, Q\}$ are free (see \cite{Spei} for freeness).
\end{itemize} 
Thus the process takes values in {\it the compressed space} $(P\mathscr{A}P, (1/\phi(P))\phi)$. The spectral distribution has the following decomposition :  
\begin{equation*}
\mu_{\lambda,\theta}(dx) = \frac{1}{2\pi \lambda \theta}\frac{\sqrt{(x_+ - x)(x-x_-)}}{x(1-x)}{\bf 1}_{[x_-,x_+]}(x)dx + a_0\delta_0(dx) + a_1\delta_1(dx)
 \end{equation*} where $\delta_y$ stands for the Dirac mass at $y$ with corresponding weight $a_y$, $y \in \{0,1\}$ and 
 \begin{equation*}
x_{\pm}  =   \left(\sqrt{\theta(1 - \lambda \theta)} \pm  \sqrt{\lambda \theta(1 - \theta)}\right)^2
\end{equation*}
Its Cauchy transform writes 
\begin{equation}\label{Cauchy}
G_{\mu_{\lambda,\theta}}(z) = \frac{(2-(1/\lambda \theta))z + (1/\lambda -1) + \sqrt{Az^2 - Bz + C}}{2z(z-1)}, \, z \in \mathbb{C} \setminus [0,1]
\end{equation}
with $A = 1/(\lambda\theta)^2$, $B= 2((1/\lambda \theta)(1+1/\lambda) -  2/\lambda)$ et $C = (1-1/\lambda)^2$. It was shown in \cite{Dem} that if $\lambda \in ]0,1], 1/\theta \geq \lambda + 1$ then the process is injective in $P\mathscr{A}P$, that is $a_0 = a_1 = 0$. Moreover, 
$\mu_{1,1/2}(dx)$ fits the Beta distribution $B(1/2,1/2)$: 
\begin{equation*}
\mu_{1,1/2}(dx) = \frac{1}{\pi \sqrt{x(1-x)}} {\bf 1}_{[0,1]}(x) dx
\end{equation*}
Recall that the Tchebycheff polynomials of the first kind are defined by 
\begin{equation*}
T_n(x) = \cos(n\arccos x), \, n \geq 0,\,  |x| \leq 1.  
\end{equation*}
and that they are orthogonal with respect to $\mu_{1,1/2}(dx)$. Their generating function is given by: 
\begin{equation*}
g(u,x) = \sum_{n \geq 0}T_n(x)u^n = \frac{1-ux}{1-2ux + u^2}, \quad |u| < 1. 
\end{equation*}
In \cite{Dem}, we proved that for $r > 0$
\begin{equation*}
g(re^t, J_t) = ((1+ re^t)P - 2e^tJ_t)((1+re^t)^2P - 4re^tJ_t)^{-1}, \quad t < -\ln r
\end{equation*}
defines a free martingale with respect to the natural filtration of $J$, say $\mathscr{J}_t$, the unit of the compressed space being the projection $P$. 
It follows that $(e^{nt}T_n(2J_t - P))_{t \geq 0},\, n \geq 1$ is a family of free martingale polynomials. Note also that 
\begin{align*}
h(re^t,J_t) &:= 2g(re^t,J_t) - P = \frac{(1- r^2e^{2t})}{(1+re^t)^2}(P -\frac{4re^t}{(1+re^t)^2}J_t)^{-1} 
\\& =  \frac{1- re^{t}}{1+re^t}(P -\frac{4re^t}{(1+re^t)^2}J_t)^{-1} 
\\& = (1- (re^t)^2)(P - 2re^t(2J_t - P) + (re^t)^2)^{-1}  
\end{align*}
is also a free martingale. Let $U_n$ denote the $n$-th Tchebycheff polynomial of the second kind defined by 
\begin{equation*}
U_n(\cos \alpha) = \frac{\sin(n+1) \alpha}{\sin \alpha}, \quad \alpha \in \re
\end{equation*}
with generating function given by 
\begin{equation*}
\sum_{n \geq 0}U_n(x)u^n = \frac{1}{1-2ux + u^2}, \quad |x| \leq 1,\,|u| < 1. 
\end{equation*}
Then, one deduces either from the above generating function or from the relation $2T_n = U_n - U_{n-2}, U_{-1}  := 0$ that $\{M_t^n := e^{nt}(U_n - U_{n-2})(2J_t - P),\, n \geq 1\}_{t \geq 0}$ is a family of free martingale polynomials. The aim of this work is to extend this claim to the range $\theta = 1/2, \lambda \in ]0,1]$. The motivation originates from \cite{Kubo2} where the author determines the family of orthogonal polynomials with respect to $\mu_{\lambda,\theta}$. Our first guess was that these will be free martingales polynomials for all $\lambda \in 
]0,1],\, \theta \leq 1/(\lambda+1)$. Yet, things turn to be more complicated: not only the range is restricted but the martingale polynomials we derive are not orthogonal with respect to 
$\mu_{\lambda,1/2}$ except for $\lambda =1$. As a matter of fact, we will on one hand derive the orthogonal polynomials corresponding to $\mu_{\lambda,1/2}$ and compute on the other hand the appropriate orthogonality measure for our martingales polynomials. The last part of the paper is devoted to a realization of the free stationary Jacobi process using the Accardi-Bozejko isomorphism (see \cite{Acca1}) as well as some comments. 

\begin{nota}
From a matrix theory point of view, the choice $\theta = 1/2$ correponds to the ultraspherical multivariate Beta distribution (see \cite{Dem}). Moreover, to our level of Knowledge, there is only one result concerning martingale polynomials for the stationary (classical) Jacobi process, which is restricted to the one dimensional case. More precisely, pick a vector $(x_1,\dots,x_d)$ belonging to the sphere $S^{d-1}, d \geq 2$ distributed according to the uniform (Haar) measure, then form the discrete process defined by 
\begin{equation*}
s_p = \sum_{i=0}^{p} x_i^2,\quad 1 \leq p \leq d-1.
\end{equation*}
It is a known that each random variable has the Beta distribution $B((d-p/2), p/2)$. It was shown in \cite{Silv} that 
\begin{equation*}
M_n^d(p) = \frac{1}{((d-p)/2)_n} \ja(2s_p - 1),
\end{equation*}
where $\ja$ denotes the $n$-th Jacobi polynomial of parameters $\alpha = (d-p)/2 -1, \beta = (p/2) -1$, is a martingale with respect to the natural filtration of the process. 
To relate this to our work, we rewrite $s_p$ in the matrix form 
\begin{equation*}
s_p = P_1U_dQ_pU_d^{\star}P_1,
\end{equation*} 
where $U_d$ is a $d \times d$ Haar unitary matrix, $P_1$ is a $d \times d$ projection with only one non vanishing coefficient $(P_1)_{11} = 1$ and $Q_p$ is a $d \times d$ projection with only $p$ non vanishing term $(Q_p)_{11} = \dots = (Q_p)_{pp} = 1$. For $d=2p$, we get the ultraspherical polynomials of parameter $\lambda = (p-1)/2$. 
\end{nota}
\section{Main result}
One has for $\lambda \in ]0,1],\, \theta = 1/2$
\begin{equation*}
x_- = \left(\frac{\sqrt{2-\lambda}}{2} - \frac{\sqrt{\lambda}}{2}\right)^2 \leq x \leq x_+ = \left(\frac{\sqrt{2-\lambda}}{2} + \frac{\sqrt{\lambda}}{2}\right)^2 \Rightarrow 
-1 \leq \frac{2x - 1}{\sqrt{\lambda(2-\lambda)}} \leq 1
\end{equation*}
and our main result is stated as follows: 
\begin{pro}
Set 
\begin{equation*}
a(\lambda) = \frac{(1-\lambda)}{\sqrt{\lambda(2-\lambda)}},\, x_{t,\lambda} = \frac{2J_t - P}{\sqrt{\lambda(2-\lambda)}}
\end{equation*}
For each $n \geq 1$, the process defined by 
\begin{equation*}
M_t^n := [U_n(x_{t,\lambda}) -  2a(\lambda)U_{n-1}(x_{t,\lambda}) - U_{n-2}(x_{t,\lambda})]\left(\frac{e^t}{\lambda(2-\lambda)}\right)^n, t \geq 0
\end{equation*} 
is a $(\mathscr{J}_t)$-free martingale.
\end{pro}

\section{Proof of the main result}
The proof consists of two parts: the first one consists in deriving a martingale function for all values of $\lambda \in ]0,1], \theta \leq 1/2 \leq 1/(\lambda + 1)$. 
In the second one, we specialize for $\theta = 1/2$ and show that this function correponds to the generating function of the polynomials stated above.\\  
{\it First step}: inspired by the above expression of $h(re^t,J_t)$, we will look for martingales of the form 
\begin{equation*}
R_t := K_t(P - Z_tJ_t)^{-1} =  K_t\sum_{n \geq 0}Z_t^nJ_t^n :=  K_tH_t
\end{equation*}  
where $K, Z$ are differentiable functions of the variable $t$ lying in some interval $[0,t_0[$ such that $0 < Z_t < 1$ for $t \in [0,t_0[$. 
The finite variation part of $dR_t$ is given by 
\begin{equation*}
FV(dR_t) = K'_tH_t dt + K_t FV(dH_t) 
\end{equation*}
Our main tool is the free stochastic calculus and more precisely the free stochastic differential equation already set for 
$J_t^n,\, n \geq 1$ (\cite{Dem}): 
\begin{equation*}
dJ_t^n = dM_t + n(\theta P - J_t)J_t^{n-1}dt + \lambda \theta \sum_{l = 1}^{n-1}l[m_{n-l}(P-J_t)J_t^{l-1}+ (m_{n-l-1} - m_{n-l})J_t^l)]dt
\end{equation*}
where $dM$ stands for the martingale part and $m_n$ is the $n$-th moment of $J_t$ in $P\mathscr{A}P$: 
\begin{equation*}
m_n := \tilde{\phi}(J_t^n) : = \frac{1}{\phi(P)} \phi(J_t^n)
\end{equation*}
The finite variation part $FV(dJ_t^n)$ of $J_t^n$ transforms to: 
\begin{align*}
FV(dJ_t^n) &= n(\theta P - J_t)J_t^{n-1}dt + \lambda \theta \left[\sum_{l = 1}^{n-1}l[m_{n-l}J_t^{l-1}+ \sum_{l=1}^{n-1} l(m_{n-l-1} - 2m_{n-l})J_t^l)\right]dt
\\& =  n(\theta P - J_t)J_t^{n-1}dt + \lambda \theta \sum_{l = 1}^{n-1}lm_{n-l}J_t^{l-1}+ \sum_{l=1}^{n} (l-1)(m_{n-l} - 2m_{n-l+1})J_t^{l-1}dt
\\& = n(\theta P - J_t)J_t^{n-1}dt + \lambda \theta \sum_{l = 1}^{n}[lm_{n-l} + (l-1)(m_{n-l} - 2m_{n-l+1})]J_t^{l-1}dt - n\lambda \theta J_t^{n-1}dt
\\& = n\theta(1-\lambda)J_t^{n-1}dt - nJ_t^n dt + \lambda \theta \sum_{l = 1}^{n}[m_{n-l} + 2(l-1)(m_{n-l} - m_{n-l+1})]J_t^{l-1}dt 
\end{align*}
Thus 
\begin{align*}
FV(dH_t) & = \sum_{n \geq 1}nZ'_tZ_t^{n-1} J_t^n dt + \sum_{n \geq 1}Z_t FV(J_t^n)
\\& = \sum_{n \geq 1}nZ'_t Z_t^{n-1}J_t^n dt -\sum_{n \geq 0}nZ_t^n J_t^n dt + \theta(1-\lambda)\sum_{n \geq 1}nZ_t^n J_t^{n-1}dt
\\& + \lambda \theta \sum_{n \geq 1}\sum_{l = 1}^n Z_t^nm_{n-l}J_t^{l-1}dt + 2\lambda \theta \sum_{n \geq 1}\sum_{l =1}^n (l-1)Z_t^n(m_{n-l} - m_{n-l+1})]J_t^{l-1}dt
\\& =  \sum_{n \geq 1}n[Z'_t Z_t^{n-1} - Z_t^n]J_t^n dt + \theta(1-\lambda)\sum_{n \geq 0}(n+1)Z_t^{n+1} J_t^{n}dt
\\& + \lambda \theta \sum_{n \geq 0}\sum_{l \geq 0} Z_t^{n+l+1}m_{n}J_t^{l}dt + 2\lambda \theta \sum_{n \geq 0}\sum_{l\geq 0} lZ_t^{n+l+1}(m_{n} - m_{n+1})]J_t^{l}dt
\end{align*}
\begin{align*}
& = [Z'_t/Z_t  - 1 + \theta(1-\lambda)Z_t]\sum_{n \geq 1}nZ_t^{n}J_t^n dt + \theta(1-\lambda)Z_t\sum_{n \geq 0}Z_t^{n} J_t^{n}dt
\\& + \lambda \theta \sum_{n \geq 0}Z_t^{n+1}m_n\sum_{l \geq 0} Z_t^{l}J_t^{l}dt + 2\lambda \theta \sum_{n \geq 0}Z_t^{n+1}(m_{n} - m_{n+1})]\sum_{l \geq 0}lZ_t^lJ_t^{l}dt
\end{align*}
Recall that the Cauchy transform of a measure on the real line is defined by  
\begin{equation*}
G_{\nu}(z) = \int_{\re}\frac{1}{z-x} \nu(dx) = \sum_{n \geq 0}\frac{1}{z^{n+1}} \int_{\re}x^n \nu(dx)
\end{equation*}
for some values of $z$ for which both the integral and the infinite sum make sense. Then, since $0 < Z < 1$ and $\mu_{\lambda,\theta}$ is supported in $[0,1]$, it is easy to see that 
\begin{equation*}
\sum_{n \geq 0}Z_t^{n+1}(m_{n} - m_{n+1}) = \left(1- \frac{1}{Z_t}\right)G_{\mu_{\lambda,\theta}}\left(\frac{1}{Z_t}\right) + 1
\end{equation*}
with $G_{\mu_{\lambda,\theta}}$ given by (\ref{Cauchy}). This gives 
\begin{equation*}
2\lambda \theta (1-z)G_{\mu_{\lambda,\theta}}(z) = \frac{(1-2\lambda \theta)z - \theta(1- \lambda) - \sqrt{z^2 - (\lambda \theta)^2Bz + (\lambda \theta)^2C}}{z},
\end{equation*}
so that 
\begin{equation*}
2\lambda \theta (1-Z_t^{-1})G_{\mu_{\lambda,\theta}}(Z_t^{-1}) + 2\lambda \theta = 1- \theta(1- \lambda)Z_t - \sqrt{1 - (\lambda \theta)^2BZ_t + (\lambda \theta)^2CZ_t},
\end{equation*}
We finally get: 
\begin{align*}
FV(dH_t) &= [Z'_t/Z_t  - \sqrt{1 - (\lambda \theta)^2BZ_t + (\lambda \theta)^2CZ_t^2}] \sum_{n \geq 1}nZ_t^{n}J_t^n dt 
\\&+ \left[\lambda \theta G_{\mu_{\lambda,\theta}}\left(\frac{1}{Z_t}\right)+ \theta(1-\lambda)Z_t\right]\sum_{n \geq 0}Z_t^{n} J_t^{n}dt
\end{align*}
In order to derive free martingales, we shall pick $Z$ such that $Z'_t = Z_t \sqrt{1 - (\lambda \theta)^2BZ_t + (\lambda \theta)^2CZ_t^2}$. This shows that $Z$ is an increasing function and one can solve the above non linear differential equation as follows: use the variables change $u = Z_t,\, t < t_0$, then integrate to get :
\begin{equation*}
\int_{[Z_0,Z_t]} \frac{du}{u \sqrt{1 - 2\theta(1+\lambda - 2\lambda \theta)u + (\theta(1-\lambda))^2u^2}} = t
\end{equation*}
\begin{nota}
Let $c_1 = 2\theta(1+\lambda - 2\lambda \theta),c_2 = \theta^2(1-\lambda)^2$. Then, the function $u \mapsto 1- c_1u + c_2u^2$ is decreasing for $u \in ]0,1[$: in fact,
\begin{align*}
2c_2u - c_1 &< 2c_2 - c_1 = 2\theta^2(1-\lambda)^2 - 2\theta(1+\lambda - 2\lambda\theta)
\\& = 2\theta[\theta(1+\lambda^2) - (1+\lambda)] \leq 2\theta\left(\frac{1+\lambda^2}{1+\lambda} - (1+\lambda)\right) = -\frac{4\lambda\theta}{1+\lambda} < 0
\end{align*}
which yields $1- c_1u + c_2u^2 > 1-c_1 + c_2 = (1- \theta(1+\lambda))^2 \geq 0$. \end{nota}
Next, use the variable change $1-vu = \sqrt{1- c_1u + c_2u^2}$. This gives 
\begin{equation*}
u = \frac{2v - c_1}{v^2 - c_2},\,du = -2\frac{v^2 +c_2 -c_1v}{(v^2 - c_2)^2}dv,\, 1-vu = -\frac{v^2 + c_2 - c_1v}{v^2 - c_2}
\end{equation*}
Moreover 
\begin{equation*}
u \mapsto v = \frac{1- \sqrt{1-c_1u + c_2u^2}}{u},\quad 0 < u < 1
\end{equation*}
is an increasing function: in fact the numerator of its derivative writes
\begin{equation*}
c_1u - 2c_2u^2 + 2(1-c_1u + c_2u^2) - 2\sqrt{1-c_1u+c_2u^2} = (2-c_1u) - 2\sqrt{1-c_1u+c_2u^2}
\end{equation*}
Since $2-c_1u > 2-c_1 = 2(1-\theta(1+\lambda)) + 4\lambda\theta^2 > 0$, our claim follows from the fact that $c_1^2-4c_2 = 16\lambda\theta^2(1-\lambda\theta)(1-2\theta) \geq 0$. 

Finally, the integral transforms to 
\begin{equation*}
\int_{[v_0,v_t]} \frac{2dv}{2v - c_2} =  \log\left|\frac{2v_t - c_1}{2v_0 - c_1}\right| = t
\end{equation*}
where $1-Z_t v_t = \sqrt{1- c_1Z_t + c_2Z_t^2}, \, 1-Z_0 v_0 = \sqrt{1- c_1Z_0 + c_2Z_0^2}$. Note also that $c_1^2 - 4c_2 \geq 0$ implies that for all $ u \in [Z_0,Z_t] \subset ]0,1[$
\begin{align*}
v - \frac{c_1}{2} &= \frac{1-\sqrt{1-c_1u + c_2u^2}}{u} - \frac{c_ 1}{2} = \frac{(1-c_1u/2) - \sqrt{1-c_1u + c_2u^2}}{u} 
\\& = \frac{(1-c_1u/2)^2 - (1-c_1u + c_2u^2)}{u((1-c_1u/2) +  \sqrt{1-c_1u + c_2u^2})}\geq 0
\end{align*}
since $1-c_1/2u \geq 1-c_1/2 \geq 0$. Thus $v \geq c_1/2 \geq \sqrt{c_2}$.   
\begin{equation*}
v_t = [(2v_0 - c_1) e^{t} + c_1]/2 \Leftrightarrow  \sqrt{1-c_1Z_t + c_2Z_t^2} = 1 - \frac{(2v_0 - c_1) e^{\pm t} + c_1}{2}Z_t
\end{equation*}
We finally get 
\begin{equation*}
Z_t = \frac{4(2v_0 - c_1)e^{\pm t}}{((2v_0 - c_1)e^{t} + c_1)^2 - 4c_2}, \quad t \leq t_0
\end{equation*}
where $t_0$ is the first time such that $Z_{t_0} = 1 \Leftrightarrow (2v_0 - c_1)e^{t_0} + c_1)^2 - 4c_2 - 4(2v_0 - c_1)e^{t_0}$. Set $r = r(\lambda,\theta) := (2v_0 - c_1)$ and 
$x_0 = e^{t_0} > 1$, then $r^2x_0^2 + 2(c_1-2)rx_0 + c_1^2 - 4c_2 = 0$. The discriminant equals to $\Delta = 16r^2(1+c_2-c_1) = 16r^2(1- \theta(1+\lambda))^2$. Thus 
\begin{equation*}
x_0 = \frac{-(c_1 - 2) - 2(1- \theta(1+\lambda))}{r} = \frac{2(1- \theta(1+\lambda)) + 4\lambda \theta^2 - 2(1-\theta(1+\lambda))}{r} = \frac{4\lambda \theta^2}{r} \geq 1
\end{equation*}  
The last inequality follows from the fact that $1- \sqrt{c_2}u \geq 1- \theta(1+\lambda) \geq 0$ and from
\begin{equation*}
r - 4\lambda \theta^2 = 2v_0 - c_1 - 4\lambda \theta^2 = 2(v_0 - \theta(1+\lambda)) = 2(v_0-\sqrt{c_2}) \leq 0.
\end{equation*}
It gives $t_0 = -\ln(r/4\lambda \theta^2)$. Note also that the denominator is well defined for all $t \leq t_0$ since $c_1^2 \geq 4c_2$ and $2v_0 - c_1 \geq 0$.\\
For the ramaining terms, we shall choose $K$ such that 
\begin{equation*}
K'_t + K_t\left[\lambda \theta G_{\mu_{\lambda,\theta}}\left(\frac{1}{Z_t}\right)+ \theta(1-\lambda)Z_t\right] = 0 
\end{equation*}
An easy computation shows that this equals to 
\begin{equation*}
K'_t + \frac{K_t}{2}\left[\theta(1-\lambda)\frac{Z_t^2}{Z_t - 1} + (1-2\theta)\frac{Z_t}{Z_t - 1} - \frac{Z_t\sqrt{1-c_1Z_t + c_2Z_t^2}}{Z_t - 1}\right] = 0
\end{equation*}
Remembering the choice of the function $Z$, this writes 
\begin{equation*}
K'_t - \frac{K_t}{2}\left[\frac{Z'_t}{Z_t - 1} - (1-2\theta)\frac{Z_t}{Z_t - 1} - \theta(1-\lambda)\frac{Z_t^2}{Z_t - 1}\right] = 0
\end{equation*}
or equialently
\begin{equation*}
K'_t - \frac{K_t}{2}\left[\frac{Z'_t}{Z_t - 1} - (1-\theta - \lambda \theta)\frac{Z_t}{Z_t - 1} - \theta(1-\lambda)Z_t\right] = 0
\end{equation*}
If $K_t \neq 0$, then 
\begin{equation*}
\log K_t = \frac{1}{2}\log (1-Z_t) - \frac{1- \theta - \lambda \theta}{2}\int \frac{Z_s}{Z_s - 1} ds - \frac{\theta(1-\lambda)}{2}\int Z_s ds + C
\end{equation*} 
If $\lambda \neq 1$, then the last term is given by 
\begin{equation*}
- \frac{\theta(1-\lambda)}{2}\int Z_s ds = \frac{\theta(1-\lambda)}{\sqrt{c_2}}\int \frac{(r/2\sqrt{c_2})e^t}{1-\displaystyle \left(\frac{re^t + c_1}{2\sqrt{c_2}}\right)^2}
= \arg \tanh \left(\frac{re^t + c_1}{2\sqrt{c_2}}\right)
\end{equation*}
where $\arg \tanh (u) = (1/2)\log((u+1)/(u-1)), |u| > 1$. 
The second term writes 
\begin{align*}
\frac{Z_t}{Z_t - 1} &=  \frac{4re^t}{4c_2 + 4re^t - (re^t + c_1)^2} =  \frac{4re^t}{4c_2 - c_1^2 + (c_1-2)^2 -  (re^t + c_1 - 2)^2}
\\& = \frac{re^t}{c_2 + 1 - c_1 - \displaystyle \left(\frac{re^t + c_1 - 2}{2}\right)^2} = \frac{1}{c_2+1-c_1} \frac{re^t}{1- \displaystyle \left(\frac{re^t + c_1 - 2}{2\sqrt{c_2+ 1 -c_1}}\right)^2}
\\& =  \frac{2}{\sqrt{c_2+1-c_1}} \frac{(r/2\sqrt{c_2+1-c_1})e^t}{1- \displaystyle \left(\frac{re^t + c_1 - 2}{2\sqrt{c_2+ 1 -c_1}}\right)^2}
\end{align*}
Observe that $2-c_1 - re^t > 2 - c_1 - re^{t_0} = 2(1- \theta(1+\lambda) \geq 0$. Thus, if $\theta(1+\lambda) \neq 1$ 
\begin{equation*}
\frac{1- \theta(1+ \lambda)}{2}\int \frac{Z_s}{Z_s - 1} ds = \arg \tanh \left(\frac{2- c_1 - re^t}{2\sqrt{c_2+ 1 -c_1}}\right)
\end{equation*}
Thus, if $\lambda \neq 1 \,(\theta \leq 1/2 < 1/(\lambda+1))$,   
\begin{align*}
K_t &= C(1-Z_t)^{1/2} \left(\frac{re^t + c_1+ 2\sqrt{c_2}}{re^t + c_1 - 2\sqrt{c_2}}\right)^{1/2}\left(\frac{2-  c_1 - 2c_3 - re^t}{2- c_1 +2c_3 - re^t}\right)^{1/2}
\end{align*}
where $c_3 := \sqrt{c_2+1-c_1} = 1- \theta(\lambda+1)$. Note that for $\lambda =1, \theta = 1/2$, $c_1 = 1, c_2 = 0, c_3 = 0$ and 
\begin{equation*}
K_t  = C\frac{1-re^t}{1+re^t}, \qquad t < t_0 = -\ln r.
\end{equation*}

{\it The case $\theta = 1/2, \lambda \neq 1$: free martingales polynomials}:
one has 
\begin{align*}
&c_1 =  1,\, c_2 = \frac{(1-\lambda)^2}{4},\, c_3 = \sqrt{c_2} = \frac{1-\lambda}{2},\, Z_t = \frac{4re^t}{(re^t + 1)^2 - (1-\lambda)^2}
\\& c_1 + 2\sqrt{c_2} = 2(1+c_3) - c_1  = 2-\lambda,\, c_1 - 2\sqrt{c_2} = 2(1-c_3) - c_1 = \lambda.\\& 
1-Z_t = \frac{(re^t - 1)^2 - (1-\lambda)^2}{(re^t + 1)^2 - (1-\lambda)^2} = \frac{(re^t + \lambda -2)(re^t - \lambda)}{(re^t +2- \lambda)(re^t + \lambda)}
\end{align*}
Thus, for $t < -\ln (r/\lambda)$, 
\begin{equation*}
K_t = C\frac{\lambda - re^t}{\lambda + re^t}
\end{equation*}
so that 
\begin{align*}
R_t &= C\frac{\lambda - re^t}{\lambda + re^t} (P - \frac{4re^t}{(re^t + 1)^2 - (1-\lambda)^2}J_t)^{-1} 
\\&= C(\lambda - re^t)(2-\lambda + re^t)(\lambda(2-\lambda) P + (re^t)^2P - 2re^t(2J_t - P))^{-1}
\\&= \frac{C(\lambda - re^t)(2-\lambda + re^t)}{\lambda(2-\lambda)} \left(P  - \frac{2re^t}{\sqrt{\lambda(2-\lambda)}} \frac{(2J_t - P)}{\sqrt{\lambda(2-\lambda)}}
+ \frac{(re^t)^2}{\lambda(2-\lambda)}P\right)^{-1}
\\& = C\left(1 - 2\frac{(1-\lambda)}{\sqrt{\lambda(2-\lambda)}}\frac{re^t}{\sqrt{\lambda(2-\lambda)}} - \frac{(re^t)^2}{\lambda(2-\lambda)}\right)
\left(P - \frac{2re^t}{\sqrt{\lambda(2-\lambda)}} \frac{(2J_t - P)}{\sqrt{\lambda(2-\lambda)}}+ \frac{(re^t)^2}{\lambda(2-\lambda)}P\right)^{-1}
 \end{align*}
is a free martingale with respect to the natural filtration $\mathscr{J}_t$.
Besides, since $\lambda \in ]0,1]$, then $\lambda \leq \sqrt{\lambda(2-\lambda)}$, hence $(re^t)/(\sqrt{\lambda(2-\lambda)}) < 1$ for all $t < -\ln(r/\lambda)$. 
Now, let us consider the following generating function
\begin{equation*}
g(u,x) = \frac{1- 2au - u^2}{1 - 2xu + u^2}, \quad 0 < a , u < 1,\,\, |x| \leq 1. 
\end{equation*}
It follows that 
\begin{equation*}
g(u,x) = U_0(x) + (U_1(x)-2a)u + \sum_{n \geq 2}[U_n(x) - 2aU_{n-1}(x) - U_{n-2}(x)]u^n
\end{equation*}
Setting 
\begin{equation*}
u_{t,\lambda} := \frac{re^t}{\sqrt{\lambda(2-\lambda)}},\quad t < t_0,
\end{equation*}
then
\begin{equation*}
R_t = C[P + (x_{t,\lambda} - 2a(\lambda)P)u_{t, \lambda} + \sum_{n \geq 2}[U_n(x_{t,\lambda}) - 2a(\lambda)U_{n-1}(x_{t,\lambda}) - U_{n-2}(x_{t,\lambda})]u_{t,\lambda}^n
\end{equation*}
Setting $U_{-1} = U_{-2} = 0$, it can be written as 
\begin{equation*}
R_t = C\sum_{n \geq 0}[U_n(x_{t,\lambda}) - 2a(\lambda)U_{n-1}(x_{t,\lambda}) - U_{n-2}(x_{t,\lambda})]u_{t,\lambda}^n 
\end{equation*}

\begin{nota} {\it The case $\lambda = 1$.}\\
$c_1 = 4\theta(1-\theta),\, c_ 2 = 0$ and $Z_t$ writes 
\begin{equation*}
Z_t = \frac{4re^t}{(re^t + 4\theta(1-\theta))^2} 
\end{equation*}
Moreover, $c_3 = \sqrt{1-c_1} = (1-2\theta),\, 2- 2c_3 - c_1 = 4\theta^2,\, 2 + 2c_3 - c_1 = 4(1-\theta)^2$. $K_t$ then writes
\begin{equation*}
K_ t = \frac{\sqrt{(re^t + 4\theta(1-\theta))^2 - 4re^t}}{re^t + 4\theta(1-\theta)} \sqrt{\frac{4\theta^2 - re^t}{4(1-\theta)^2 - re^t}}
\end{equation*}
\end{nota}

\section{one-parameter measures family and Orthogonal polynomials}
Let $\mu$ be a measure on the real line which is not supported by a finite set. Assume that $\mu$ has finite moments of all orders. Applying the Gram-Schmidt orthogonolization method to the basis $(1,x,x^2,\dots)$, there exist a unique family of monic orthogonal polynomials with respect to $\mu$, say $(P_n)_{n \geq 0}$. These polynomials satisfy the three-terms recurrence relation 
\begin{equation*}
(x- \alpha_n)P_n(x) = P_{n+1}(x) + \omega_nP_{n-1}(x),\quad n \geq 0,  P_{-1} := 0.
\end{equation*}
where $\alpha_n \in \re,\, w_n > 0$. $(\alpha_n,\omega_n)_{n \geq 0}$ are called the Jacobi-Szeg\"o parameters of $\mu$. It is known that $\mu$ is symmetric if and only if 
$\alpha_n = 0,\, n \geq 0$. 
Another way to derive the family $(P_n)_n$ is the multiplicative renormalization method (\cite{Asai1},\cite{Asai1b},\cite{Asai2}, \cite{Asai3}) that we shall recall here :  a nice function $(u,x) \mapsto \psi(u,x)$ is a generating function for the measure $\mu$ if $\psi$ has the expansion 
\begin{equation*}
\psi(u,x) = \sum_{n \geq 0}c_nP_n(x)u^n, \qquad c_n \in \re 
\end{equation*}
where $P_n$ are orthogonal with respect to $\mu$. Of course, there is more than one generating function corresponding to a given measure and in order to claim whether a function 
is a generating function or not, authors in \cite{Asai1} provided a necessary and sufficient condition. 
For a particular form of $\psi$ which fits our need, their result is formulated as follows: 
\begin{teo}\label{T1}
Define  
\begin{equation*}
\theta(u) := \int_{\re}\frac{1}{1-ux}\mu(dx), \quad \theta(u,v) := \int_{\re}\frac{1}{(1-ux)(1-vx)}\mu(dx).
\end{equation*}
Let $\rho$ analytic around $0$ such that $\rho(0) = 0$ and $\rho'(0) \neq 0$. Then 
\begin{equation}\label{generating}
\psi(u,x) := \frac{(1-\rho(u)x)^{-1}}{\theta(\rho(u))}
\end{equation}
is a generating function for $\mu$ if and only if 
\begin{equation*}
\Theta_{\rho}(u,v)  := \frac{\theta(\rho(u),\rho(v))}{\theta(\rho(u))\theta(\rho(v))}
\end{equation*}
is a function of $uv$.
\end{teo}   
We will apply this result to the measures family $\nu_{\lambda}, \lambda \in ]0,1]$ which is the image of 
\begin{equation*}
\mu_{\lambda,1/2} = \frac{1}{\pi \lambda}\frac{\sqrt{(x_+ - x)(x-x_-)}}{x(1-x)}{\bf 1}_{[x_-,x_+]}(x)dx,\quad x_{\pm} = \frac{(\sqrt{\lambda} \pm \sqrt{2-\lambda})^2}{4}
\end{equation*} 
by the map 
\begin{equation*}
x \mapsto \frac{2x - 1}{\sqrt{\lambda(2-\lambda)}}
\end{equation*}
Then,
\begin{equation*}
\nu_{\lambda}(dx) = \frac{(2-\lambda)}{\pi}\frac{\sqrt{1-x^2}}{1 - \lambda(2-\lambda)x^2}{\bf 1}_{[-1,1]}(x)dx
\end{equation*}
Our scheme is the almost the same used in \cite{Kubo1} except the computation of $\theta(u)$ which follows easily from $G_{\mu_{\lambda,1/2}}$. More precisely, authors considered the one-parameter measures family
\begin{equation*}
\mu_{a}(dx) = \frac{a\sqrt{1-x^2}}{a^2 + (1-2a)x^2} {\bf 1}_{]-1,1[} dx,\quad a > 0.
\end{equation*}
It is forward that $\mu_{1/(2-\lambda)} = \nu_{\lambda}$ almost everywhere for $0 < \lambda \leq 1 \Leftrightarrow 1/2 < a \leq 1$.
\begin{pro}
\begin{equation*}
\theta(u) = \theta_{\lambda}(u) = \frac{2-\lambda}{1- \lambda + \sqrt{1-u^2}}, \quad |u| < 1
\end{equation*}
\end{pro}
Using 
\begin{equation*}
\frac{1}{(1-ux)(1-vx)} = \frac{1}{u - v}\left(\frac{u}{1-ux} - \frac{v}{1-vx}\right)
\end{equation*}
it follows that $\theta(u,v)  = (u\theta(u) - v\theta(v))/(u-v)$ from which we deduce
\begin{cor}
\begin{equation*}
\theta(u,v) = \theta_{\lambda}(u,v) = \frac{1}{2-\lambda}\left[1-\lambda + \frac{u+v}{u\sqrt{1-v^2} + v\sqrt{1-u^2}}\right]
\end{equation*}
\end{cor}
{\it Proof:} from the definition of $\nu_{\lambda}$, one writes for $ 0 <u < \lambda(2-\lambda) \leq 1$:
\begin{equation*}
\int_{\re}\frac{1}{1-ux}\nu_{\lambda}(dx) = \int_{\re}\frac{1}{1-u\displaystyle \frac{2x-1}{\sqrt{\lambda(2-\lambda)}}} \mu_{\lambda,1/2}(dx) 
= \frac{\sqrt{\lambda(2-\lambda)}}{2u}G_{\mu_{\lambda,1/2}}\left(\frac{\sqrt{\lambda(2-\lambda)}+u}{2u}\right)
\end{equation*}
The result follows from 
\begin{equation*}
G_{\mu_{\lambda,1/2}}(z) = \frac{(1-\lambda)(2z-1) - \sqrt{4z^2 - 4z + (1-\lambda)^2}}{2\lambda z(1-z)},\quad z \in \mathbb{C} \setminus [0,1] \qquad \blacksquare
\end{equation*}

Let $\rho(u) = 2u/(1+u^2)$, then 
\begin{equation*}
\frac{\rho(u) + \rho(v)}{\rho(u)\sqrt{1-\rho^2(v)} + \rho(v)\sqrt{1-\rho^2(u)}} = \frac{1+ uv}{1-uv}
\end{equation*} 
so that Theorem \ref{T1} applies and claims that 
\begin{equation*}
\psi_{\lambda}(u,x) = \frac{1 - \lambda/(2-\lambda) u^2}{1 - 2ux + u^2}
\end{equation*}
is a generating function for $\nu_{\lambda}$ corresponding to the polynomials 
\begin{eqnarray*}
Q_n^{\lambda}(x) &=& U_n(x) - \frac{\lambda}{2-\lambda}U_{n-2}(x),, \quad n \geq 0,\, U_{-1} = U_{-2} := 0.
\end{eqnarray*}
Using the recurrence relation 
\begin{equation}\label{Tcheb}
2xU_n(x) = U_{n+1}(x) + U_{n-1}(x),\quad U_{-1} := 0,
\end{equation}
These polynomials satisfy 
\begin{eqnarray*}
2xQ_0^{\lambda}(x) &=& Q_1^{\lambda}(x) \\
2xQ_1^{\lambda}(x) &=& Q_2^{\lambda}(x) + \left(1+ \frac{\lambda}{2-\lambda}\right) Q_0^{\lambda}(x)\\
2xQ_n^{\lambda}(x) &=& Q_{n+1}^{\lambda}(x) + Q_{n-1}^{\lambda}(x),\, n \geq 2.
\end{eqnarray*}
Setting $Q_{-1}^{\lambda}:= 0$ and since the coefficient of the leading power in $Q_n^{\lambda}(x)$ is $2^n$, then one deduces that the Jacobi-Szeg\"o parameters are given by : 
$\alpha_n = 0,\, n \geq 0,\,w_1 = 1/(2(2-\lambda)),\, w_n = 1/4,\, n \geq 2$.  

\begin{nota}
In \cite{Kubo2}, authors characterize the absolutely continuous measures for which the multiplicative renormalization method is applicable with the generating function given by 
(\ref{generating}). They derived a two-parameters densities family written as  
\begin{equation*}
f(x) = \frac{c\sqrt{1-x^2}}{\pi[b^2 + c^2 -2b(1-c)x + (1-2c)x^2]} {\bf 1}_{[-1,1]}(x),\quad |b| < 1-c, \, 0 < c \leq 1.
\end{equation*}
These densities fit the image of absolutely continuous part of $\mu_{\lambda,\theta}$ by the map 
\begin{equation*}
u = \frac{2x - s}{d} \in [-1,1]
\end{equation*}
with $d = d(\lambda,\theta) = x_+ - x_- = 4\theta\sqrt{\lambda(1-\theta)(1-\lambda \theta)} ,\,, s = s(\lambda,\theta) = x_+ + x_- = 2\theta(1+\lambda - 2\lambda \theta)$. 
One gets 
\begin{equation*}
\nu_{\lambda,\theta}(dx) = \frac{d^2}{2\pi \lambda \theta} \frac{\sqrt{1-x^2}}{s(2-s)  + 2d(1-s) x - d^2x^2}dx
\end{equation*}
which provides the following relations 
\begin{equation}\label{relations}
c = \frac{1}{2(1-\lambda \theta)}, \quad b = \sqrt{\frac{\lambda}{(1-\theta)(1-\lambda \theta)}}(2\theta-1) 
\end{equation}
As a result, one can derive the correponding orthogonal polynomials for $\lambda \in ]0,1], \theta \leq 1/(\lambda+1)$ from the generating function (\cite{Kubo2}):
\begin{equation}\label{Kuo}
\phi(u,x) = \frac{1 - 2bu +(1-2c)u^2}{1 - 2ux + u^2}.
\end{equation} 
\end{nota}

\section{more orthogonal polynomials}
Consider the polynomials $P_n^{\lambda}$ defined by
\begin{equation*}
P_n^{\lambda}(x) = U_n(x) - 2a(\lambda)U_{n-1}(x) - U_{n-2}(x), \quad U_{-1} = U_{-2} := 0
\end{equation*}
with generating function 
\begin{equation*}
g(u,x) = \frac{1-2a(\lambda)u -u^2}{1-2xu + u^2}, \quad a(\lambda) = \frac{1-\lambda}{\lambda(2-\lambda)},\, 0 < u < 1.
\end{equation*}
The $P_n^{\lambda}$'s appear in \cite{Al Salam} as a limiting case of the $q$-Pollaczek polynomials. The coefficient of the highest monomial is equal to $2^n$. Using (\ref{Tcheb}), one deduces that 
\begin{eqnarray*}
2[x - a(\lambda)]P_0^{\lambda}(x)&=& P_1^{\lambda}(x)  \\
2xP_1^{\lambda}(x) &=& P_2^{\lambda}(x) + 2 P_0^{\lambda}(x)\\
2xP_n^{\lambda}(x) &=& P_{n+1}^{\lambda}(x) + P_{n-1}^{\lambda}(x),\quad n \geq 2.
\end{eqnarray*}
Thus the Jacobi-Szeg\"o parameters are given by $\alpha_0 = a(\lambda)$ and $\alpha_n = 0$ for all $n \geq 1$ and $\omega_1 = 1/2, \,\omega_n  = 1/4,\, n \geq 2$ 
($P_{-1}^{\lambda} = 0$).\\ 
One can use Theorem \ref{T1} to determine the probability measure, $\xi_{\lambda}$, with respect to which the $P_n^{\lambda}$s are orthogonal. Since $\alpha_0 \neq 0$, then $\xi_{\lambda}$ is not symmetric. Indeed, keeping the same function $\rho$ previously defined, then the function $\theta$ must be equal to 
\begin{equation*}
\theta(\rho(u)) = \frac{1+u^2}{1-2a(\lambda) - u^2}
\end{equation*}
so that 
\begin{equation*}
\theta(u) = \frac{1}{\sqrt{1-u^2} - a(\lambda)u}
\end{equation*}
From the definition of $\theta$, one deduces that 
\begin{equation*}
G_{\xi_{\lambda}}(u) := \int_{\re}\frac{1}{u-x}\xi_{\lambda}(dx) = \frac{1}{u}\theta\left(\frac{1}{u}\right) = \frac{\sqrt{u^2-1} + a(\lambda)}{u^2 - (1+a^2(\lambda))}
\end{equation*}
for $|u| > 1,\, u \neq \pm\sqrt{1+a(\lambda)^2}$. 
Thus, $\xi_{\lambda}$ has two atoms $a_{\pm}$ at $\pm \sqrt{a^2(\lambda) +1}$ and an absolutely continuous part given by 
\begin{equation*}
a_{\pm} = -\lim_{y \rightarrow 0^+} y\Im G_{\xi_{\lambda}}(\pm \sqrt{a^2(\lambda)+ 1} + iy),\quad g(x) = -\frac{1}{\pi}\lim_{y \rightarrow 0^+} \Im G_{\xi_{\lambda}}(x+iy)
\end{equation*}
Using that the Cauchy transform maps $\mathbb{C}^+$ to $\mathbb{C}^-$, one finally gets
\begin{equation*}
\xi_{\lambda}(dx) = \frac{a(\lambda)}{\sqrt{a^2(\lambda) + 1}} \delta_{\sqrt{a^2(\lambda) + 1}}(dx)  + \frac{1}{\pi}\frac{\sqrt{1-x^2}}{a^2(\lambda) + 1 -x^2}{\bf 1}_{|x| < 1} dx
\end{equation*}
\begin{nota}
To see that this defines a probability measure for $\lambda \neq 1$, it suffices to write
\begin{align*}
\frac{1}{\pi}\int_{-1}^1\frac{\sqrt{1-x^2}}{a^2(\lambda) + 1 -x^2}dx &= \frac{1}{\pi}\int_{0}^1\frac{\sqrt{1-x}}{\sqrt{x}(a^2(\lambda) + 1 -x)}dx
\\& = \frac{1}{2(a^2(\lambda) + 1)}{}_2F_1\left(1,\frac{1}{2},2;\frac{1}{a^2(\lambda)+1}\right) 
\end{align*}
where ${}_2F_1$ denotes the Gauss hypergeometric function given by 
\begin{equation*}
{}_2F_1(e,b,c;z) = \frac{\g(c)}{\g(b)\g(c-b)} \int_0^1 x^{b-1}(1-x)^{c- b -1}(1 - zx)^{-e} dx,\,\, \Re(b) \wedge \Re(c-b) > 0
\end{equation*}
for $|u| < 1$. Then, one uses the identity
\begin{equation*}
{}_2F_1(1,b,2;z)  = \frac{1- (1-z)^{1-b}}{(1-b)z} 
\end{equation*}
to get 
\begin{equation*}
\frac{1}{\pi}\int_{-1}^1\frac{\sqrt{1-x^2}}{a^2(\lambda) + 1 -x^2}dx = 1 - \frac{a(\lambda)}{\sqrt{a^2(\lambda)+1}}
\end{equation*}
\end{nota}

\section{One mode Interacting Fock space}
In the sequel, we give a realization of $\nu_{\lambda,\theta}$, image of the spectral measure $\mu_{\lambda,\theta}$ for $\lambda \in ]0,1], \theta \leq 1/(\lambda + 1)$ so that the support is $[-1,1]$. In the quantum scope, it is known as the quantum decomposition of $\nu_{\lambda,\theta}$. We only need the Jacobi-Szeg\"o parameters in order to apply Accardi-Bozejko Theorem (\cite{Acca1}). We first write down from the generating function (\ref{Kuo}) the orthogonal polynomials (see \cite{Kubo2}) corresponding to $\nu_{\lambda,\theta}$: 
\begin{equation*}
Q_n^{\lambda,\theta} = U_n - 2bU_{n-1} + (1-2c)U_{n-2},\quad U_{-1} = U_{-2} = 0,
\end{equation*}
where $b = b(\lambda,\theta), c = c(\lambda,\theta)$ are given by (\ref{relations}). It follows that $\alpha_0 = b, \alpha_n = 0$ for $n \geq 1$ and $\omega_1 = c/2, \omega_n = 1/4$ for 
$n \geq 1$. In order to use Accardi-Bozejko Theorem (\cite{Acca1}), we shall introduce the so-called \emph{one-mode interacting Fock space}: let $\mathcal{H}$ be a one dimensional separable complex Hilbert space $\sim \mathbb{C}$. Then the $n$-th tensor product $\mathcal{H}^{\otimes n}$ is one dimensional: indeed $z_1\otimes \dots \otimes z_n = 
(z_1\dots z_n) 1 \otimes \dots \otimes 1 \in \mathbb{C} \Phi_n$. The one-mode interacting Fock space associated to $\nu_{\lambda,\theta}$ is defined as $\Gamma(\mathbb{C} \Phi_n,(\lambda_n))$ as the infinite orhogonal sum of $\mathbb{C} \Phi_n$ equipped with the weighted scalar product 
\begin{equation*}
(z_1\Phi_n,z_2\Phi_n) := \lambda_n\overline{z_1}z_2,\quad z_1,z_2 \in \mathbb{C},
\end{equation*}        
where $\lambda_n = \omega_1\dots \omega_n$. Then $\nu_{\lambda,\theta}$ is the vacuum distribution (in the vacuum state $\Phi_1$) of any extension of the operator 
$a^+ + a + \alpha_N$ where
\begin{eqnarray*}
a^+ \Phi_n &=& \Phi_{n+1} \quad \textrm{(creation operator)} \\
a \Phi_{n+1} & = & \omega_{n+1}\Phi_n = \frac{\lambda_{n+1}}{\lambda_n}\Phi_n,\,a\Phi_1 = 0,\quad \textrm{(annihilation operator)} \\
N\Phi_n &=& n \Phi_n \qquad \textrm{(Number operator)}, \quad aa^+\Phi_n = \frac{\lambda_{n+1}}{\lambda_n} \Phi_n,  
\end{eqnarray*}  
and $\alpha_N$ is defined by the spectral Theorem, that is $\alpha_N\Phi_n = \alpha_n \Phi_n$.
\begin{nota}
The concept of one mode interacting Fock space (IFS) is purely algebraic as the reader can see from \cite{Acca1} and is fully characterized by both the commutation relations between creation and annihilation operators and $a\Phi_1 = 0$. The most important feature of Accardi-Bozejko Theorem is illustrated in the \emph{canonical} isomorphism between one mode IFS and the $L^2$-space of a given measure $\mu$ of all order moments. It is noteworthy that only the $\omega_n$s are involved in the commutation relations (thus in both one mode IFS and $L^2(\mu)$) while the $\alpha_n$s reflect only the symmetry of $\mu$.
\end{nota}         

{\bf Acknowledgments} :  The author want to thank Professors K. Dykema and M. Anshelevich as well as the organization team for their financial support to attend the concentration week on probability and analysis at Texas A$\&$M university where the author started this work and hospitality. A special thank for Professor P. Graczyk who invites the author to the 28-th conference on quantum probability and related topics, and to Professor R. Quezada Batalla for the financial support to attend the conference and talk about this work. The author is grateful to Professor L. Accardi for detailed explanations about one mode IFS and to Professor H. H. Kuo for useful remarks on the manuscript.


\begin{thebibliography}{99}
\bibitem{Acca1}\emph{L. Accardi, M. Bozejko}. Interacting Fock space and Gausssianization of probability measures. {\it Infin. Dimens. Anal. Quantum Probab. Relat. Top.1.} {\bf 4}. 1998, 663-670.
\bibitem{Al Salam} \emph{W. A. Al-Salam, T. S. Chihara}. $q$-Pollaczek polynomials and a conjecture of Andrews and Askey. {\it SIAM J. Math. Anal.} {\bf 18} (1987), no. 1, 228--242.
\bibitem{Asai1}\emph{N. Asai, I. Kubo, H. H. Kuo.} Multiplicative renormalization and generating function I. {\it Taiwanese J. Math. 7}. {\bf 1}. 2003, 89-101.
\bibitem{Asai1b}\emph{N. Asai, I. Kubo, H. H. Kuo.} Renormalization, orthogonalization and generating function. {\it Quantum Information V}. 2006, 49-55. World Sci. Publ. Hackensack, NJ. 
\bibitem{Asai2}\emph{N. Asai, I. Kubo, H. H. Kuo.}Multiplicative renormalization and generating function II. {\it Taiwanese J. Math. 8}. {\bf 4}. 2004, 593-628.
\bibitem{Asai3}\emph{N. Asai, I. Kubo, H. H. Kuo.} Generating function method for orthogonal polynomials and Jacobi-Szeg\"o parameters. {\it Probab. Math. Statist. 23}. {\bf 2}, 2003, 273-291. Acta Univ. Wratislav. No 2593. 
\bibitem{Bia}\emph{P. Biane}. Free Brownian motion, free stochastic calculus and random matrices. {\it Fie. Inst. Comm.}  {\bf 12}, Amer. Math. Soc. Providence, RI,  1997, 1-19.
\bibitem{Dem}\emph{N. Demni}. Free Jacobi process. To appear in {\it J. Theo. Proba.}
\bibitem{Kubo1}\emph{I. Kubo, H. H. Kuo, S. Namli.} Interpolation of Chebyshev polynomials and interacting Fock spaces. {\it Infin. Dimens. Anal. Quantum Probab. Relat. Top. 9.} {\bf 3}. 2006, 361-371. 
\bibitem{Kubo2}\emph{I. Kubo, H. H. Kuo, S. Namli}. The class of measures applicable to the renormalization method for $(1-x)^{-1}$. Preprint 2007.
\bibitem{Silv}\emph{M. L. Silverstein}. Orthogonal polynomial martingales on spheres. {\it S\'em. Probab. XX}. 1984/85, 419-422. Lecture notes in Math., {\bf 1204}, Springer, Berlin, 1986.  
\bibitem{Spei}\emph{R. Speicher}. Combinatorics of Free Probability Theory. {\it Lectures. I. H. P. Paris}. 1999. 
\end{thebibliography}
\end{document}